\documentclass[11pt]{article}
\usepackage{amssymb}
\usepackage{latexsym}
\newtheorem{thm}{Theorem}[section]
\newtheorem{lem}[thm]{Lemma}
\newtheorem{cor}[thm]{Corollary}
\newtheorem{pro}[thm]{Proposition}
\newtheorem{ex}[thm]{Example}

\newtheorem{defi}[thm]{Definition}

\newcommand{\gm }{\Gamma }
\newcommand{\lon }{\longrightarrow }
\newcommand{\be }{\begin{eqnarray*}}
\newcommand{\ee }{\end{eqnarray*}}

\newcommand{\per }{\backl }

\newcommand{\poidd }[2]{#1\gpd #2}
\newcommand{\poiddd }[3]{ (#1\gpd #2, \alpha_{#3}, \beta_{#3})}

\input amssym.def
\input amssym

\newcommand{\pf}{\noindent{\bf Proof.}\ }
\newcommand{\qed}{\begin{flushright} $\Box$\ \ \ \ \ \
                  \end{flushright}}

\newcommand{\reals}{{\mathbb R}}

\newcommand{\integers}{{\mathbb Z}}

\newcommand{\half}{\frac{1}{2}}

\newcommand{\backl}{\mathbin{\vrule width1.5ex height.4pt\vrule height1.5ex}}

\newcommand{\calf}{{\cal F}}
\newcommand{\calg}{{\cal G}}

\newcommand{\calx}{{\cal X}}

\newcommand{\smalcirc}{\mbox{\tiny{$\circ $}}}

\def\description label#1{\hfil\bf[#1]\hfil}
\newcommand{\diff}[1]{\frac{d}{d#1}}
\newcommand{\pttt}[2]{\frac{\partial #1}{ \partial #2}}
\newcommand{\ttt}{t_{0}}
\newcommand{\gmex}{U_{ex}(\gm )}
\newcommand{\gmm}{U_{0}(\gm )}
\newcommand{\eex}{U_{0}(E)}
\newcommand{\frakf}{{\cal F}}

\def\sdp{\mathbin{\hbox{$\mapstochar\kern-.3333em\times$}}}
\def\pds{\mathbin{\hbox{$\times\kern-.55em\mapstochar\,$}}}

\newcommand{\wed}{\mathbin{\lower1.5pt\hbox{$\scriptstyle{\wedge}$}}}

\let\Tilde=\widetilde
\let\Bar=\overline

\def\chigh{{\raise1.5pt\hbox{$\chi$}}}
\let\phi=\varphi

\def\til0{\Tilde{0}}

\def\dminus{\raise2pt\hbox{\vrule height1pt width 2ex}\hskip3pt}

\def\pback#1{\mathbin{{{\lower1.2ex\hbox{$\times$}}\atop #1}}}

\def\vlra{\hbox{$\,-\!\!\!-\!\!\!-\!\!\!-\!\!\!-\!\!\!
-\!\!\!-\!\!\!-\!\!\!-\!\!\!-\!\!\!\longrightarrow\,$}}

\def\gpd{\,\lower1pt\hbox{$\longrightarrow$}\hskip-.24in\raise2pt
             \hbox{$\longrightarrow$}\,}

\def\lgpd{\,\lower1pt\hbox{$\vlra$}\hskip-1.02in\raise2pt\hbox{$\vlra$}\,}

\def\llgpd{\,\lower1pt\hbox{$\vvlra$}\hskip-1.3in\raise2pt\hbox{$\vvlra$}\,}

\begin{document}

\title{{\bf Flux homomorphism on symplectic groupoids}
\thanks{1991 {\em Mathematics
Subject Classification.} Primary 58F05. Secondary 22E65, 22A22, 58B25, 58H05.}}

\author{ PING XU \thanks{Research partially supported by NSF
	grants  DMS92-03398 and DMS95-04913, and a NSF postdoctoral fellowship.}\\ 
Department of Mathematics\\The  Pennsylvania State University \\
University Park, PA 16802, USA\\
	{\sf email: ping@math.psu.edu }}

\date{}

\maketitle

\begin{abstract}
For any  Poisson manifold $P$, the Poisson bracket on $C^{\infty}(P)$
extends to   a  Lie bracket on the space $\Omega^{1}(P)$
of   all differential one-forms, under which the space $Z^{1}(P)$
of closed one-forms and the space $B^{1}(P)$ of
exact one-forms are Lie subalgebras. These Lie algebras
are related by the exact sequence:
$$0\lon \reals \lon
 C^{\infty}(P)\stackrel{d}{\lon } Z^{1}(P)\stackrel{f}{\lon }
 H^{1}(P, \reals )\lon 0, $$
where $ H^{1}(P, \reals )$ is considered as a trivial Lie
algebra, and $f$ is the map sending each closed one-form to its
cohomology class.
The  goal of the present paper is
 to lift this  exact sequence to
the group level for  compact Poisson manifolds under certain
integrability condition.
In particular, we  will  give a  geometric description
of a Lie group integrating the underlying Poisson algebra
$C^{\infty}(P) $.
The  group homomorphism obtained by lifting $f$
is called the  flux homomorphism for symplectic
groupoids,
which can be considered as a generalization,
  in the context of Poisson manifolds, of
 the usual  flux homomorphism   of the symplectomorphism groups
of symplectic manifolds
introduced by Calabi.

\end{abstract}
\section{Introduction}
Poisson manifolds are a  natural generalization
of symplectic manifolds. In fact, one
can think of  a Poisson manifold as a union
of symplectic manifolds of varied dimensions fitting together
in a  certain smooth way. A Poisson manifold, by
definition, is  a manifold $P$
on which the space of functions $C^{\infty}(P)$
carries  a Poisson bracket: a Lie bracket satisfying
the Lebniz rule $\{f, gh\}=\{f, g\}h+\{f, h\}g$ \cite{Lichnerowicz:1977}
\cite{Weinstein:1983}.
The underlying Lie algebra on $C^{\infty}(P)$  is then
called the Poisson algebra.
To each Poisson manifold $P$, there   is 
an antisymmetric contravariant 2-tensor $\pi$, called 
the Poisson tensor, 
such that the Poisson bracket is given by
$\{f, g\}=\pi (df , dg )$ (see  \cite{Vaisman:book} \cite{Weinstein:1983}
for more details on this subject).

 Associated to a  given  Poisson manifold $P$, there are various
 infinite dimensional Lie algebras, especially
 the Poisson algebra $C^{\infty}(P)$.
It is very  natural to ask if there exist  Lie groups
integrating these Lie algebras, and if yes,  how these
Lie groups can be described geometrically and how  they are
related.
  If $P$ is symplectic, closely related to the Poisson algebra $C^{\infty}(P)$ 
are the Lie algebra of symplectic vector fields $\calx_{s}(P)$ and
that  of all hamiltonian vector fields 
$\calx_{h}(P)$. Clearly,  $C^{\infty}(P)$ is 
a one-dimensional  central extension 
of $\calx_{h}(P)$. Globally,  the group 
$Diff(P, \omega )$ of all symplectomorphisms  has $\calx_{s}(P)$ as its
 Lie algebra,
and its  subgroup $Ham(P, \omega )$
consisting of hamiltonian symplectomorphisms
can be considered as a  Lie group integrating  $\calx_{h}(P)$
 \cite{Banyaga}  \cite{Calabi} \cite{McduffS}. 
Then a  Lie group integrating $C^{\infty}(P)$   will be 
  a one-dimensional  central extension
of  $Ham(P, \omega )$, which in fact  can be  interpreted  as the group
of all contact diffeomorphisms of the  prequantum bundle 
over $P$ when $P$ is prequantizable \cite{Brylinski:book}
\cite{Weinstein:mathZ}. 
The flux homomorphism  $F: \widetilde{Diff_{0}(P, \omega )}\lon
H^{1}(P, \reals )$ was first introduced by Calabi \cite{Calabi}
as an invariant   characterizing the hamiltonian symplectomorphisms,
 and was intensively studied by Banyaga \cite{Banyaga} and McDuff \cite{Mcduff} \cite{McduffS}.
Here $Diff_{0}(P, \omega )$ denotes the identity component of
$Diff(P, \omega )$ and $ \widetilde{Diff_{0}(P, \omega )}$
refers to  its universal covering. 
Since we will use the  flux homomorphism repeatedly in the
paper, we recall
its definition below.

  Any isotopy $\phi_{t}$ of symplectomorphisms
in  $Diff_{0}(P, \omega )$ corresponds to a time-dependent
vector field $Z_{t}$ defined by:
\begin{equation}
\label{eq:Z}
\diff{t} \phi_{t}=Z_{t}\smalcirc \phi_{t}.
\end{equation}
The flux homomorphism is the map
$F: \widetilde{Diff_{0}(P, \omega )}\lon H^{1} (P, \reals )$
defined  by
\begin{equation}
\label{eq:flux}
\{\phi_{t}\}\lon \int_{0}^{1}[Z_{t}\per \omega ] dt,
\end{equation}
where $\omega$ is the symplectic structure on $P$.
Intrinsically, the flux homomorphism 
  can be understood as the  lifting, to the group level, 
 of the Lie algebra homomorphism $f: \calx_{s}(P)\lon H^{1}(P, \reals )$
 given by:
$$ f: Z\lon [Z\per \omega ]. $$
While the kernel of $f$ is  the Lie subalgebra  $\calx_{h}(P)$,
 the kernel of $F$ should characterize  the 
group of  hamiltonian symplectomorphisms \cite{Banyaga} \cite{McduffS}.

When $P$ is a general Poisson manifold, 
the Lie algebra $C^{\infty}(P)$ could be very
 far from either the Lie algebra
of Poisson vector fields or the Lie algebra
of hamiltonian vector fields.
   For example, when $P$ is a zero Poisson manifold,
which is  an extreme case, the Lie algebra
of Poisson vector fields consists of  all vector
fields but the Lie algebra
of hamiltonian vector fields   contains only  the  zero element.
Both of these Lie algebras  are quite different from the Poisson
algebra $C^{\infty}(P)$.
It turns out that, for Poisson manifolds,  it is  more interesting
 to consider  various  Lie algebras
  of differential one-forms, instead of  those of vector fields.
In the case of symplectic manifolds, the space of 
one-forms is  isomorphic  to that of vector fields as
Lie algebras. Therefore, everything will reduce 
to the previous case.

Recall that the Poisson  tensor $\pi$
on a Poisson manifold $P$ naturally
induces a Lie algebra structure on
the space of all one-forms $\Omega^{1}(P)$, which
is given  by:
\begin{equation}
\label{eq:poisson-algebroid}
[\omega, \theta ]=L_{X_{\omega}}\theta - L_{X_{\theta}}\omega -
d[\pi (\omega, \theta)],
\end{equation}
where $X_{\omega}=\pi^{\#}\omega$ and $X_{\theta}=\pi^{\#}\theta$.
Here $\pi^{\#}$ is the bundle map from $T^{*}P$ to $TP$
defined by $<\pi^{\#}(\omega ), \theta > =\pi (\omega ,  \theta )$,
for any $\omega , \ \theta \in T^{*}P$.
In fact, $T^*P$ together with  this bracket and the
anchor map $\rho=\pi^{\#} :T^{*}P \lon TP$ becomes
a Lie algebroid \cite{CDW} \cite{Weinstein:1987}. 
It is simple to see
that both the  space  $Z^{1}(P)$ of   closed one-forms and 
the space  $B^{1}(P) $ of exact one-forms are closed
under this bracket, and therefore   they
 are Lie subalgebras.
On the Lie algebra level, one has the following 
exact sequence of Lie algebra homomorphisms:

\begin{equation}
\label{eq:exact}
0\lon \reals \lon
 C^{\infty}(P)\stackrel{d}{\lon } Z^{1}(P)\stackrel{f}{\lon }
 H^{1}(P, \reals )\lon 0,
\end{equation}
where $ H^{1}(P, \reals )$ is considered as a trivial Lie
algebra, and $f$ is the map sending each closed one-form to its
cohomology class. Our goal of the  paper, roughly
speaking, is to
lift this  exact sequence to an exact sequence on the
level of groups, and to
 analyze  the  relation between their 
groups. The essential object for the construction
is the so called {\em symplectic groupoid}. 
Since the notion of
symplectic groupoids appears  almost everywhere in the paper,
let us briefly  recall it below.
A symplectic groupoid is a Lie groupoid $\poidd{\gm}{P}$
 (always assumed to be Hausdorff in our case) with  a compatible symplectic
structure in the sense that the graph of
multiplication: $\Lambda=\{(x, y, z)|z=xy, \beta (x)=\alpha (y)\}$
is a lagrangian submanifold of $\gm \times \gm \times \Bar{\gm}$ (the bar
means taking the opposite symplectic structure).
It is well known  that the base manifold  $P$ of
a symplectic groupoid   naturally carries a  Poisson structure. 
If a Poisson manifold $P$ can be realized as 
the base Poisson manifold of
 some symplectic groupoid, it is   called integrable.
Not every Poisson manifold is  integrable.
However, a fairly large class of Poisson manifolds, for example, Lie-Poisson
structures,  are known to be integrable.
On the other hand, it was a theorem of Karasev and Weinstein that
every Poisson manifold admits a local
symplectic groupoid (we refer the reader to \cite{CDW} 
\cite{Vaisman:book} \cite{Weinstein:1987} for more details on the subject).

The relevance of of the symplectic groupoid for this
kind of construction  was also
noted some time ago by several people including Dazord, Karasev,
Weinstein. The reason for such a relevance  is quite simple.
 The Lie algebra of one-forms is very
 closely related to the
Lie algebroid $T^*P$, and the  global object
of  the latter is exactly the
symplectic groupoid $\gm$ over $P$.
So  it is not surprising that   Lie groups  of
$Z^{1}(P)$ and $B^{1}(P)$  can
  be described in terms of the groupoid $\gm$.
In fact, it is known  that  the group $U(\gm )$ of lagrangian bisections
of $\gm$ (see Section 2 for the definition of lagrangian
bisections)
is   a group integrating $Z^{1}(P)$.

On the group level,  the morphism $f$ should correspond to a group
homomorphism $\frakf :\widetilde{\gmm}\lon H^{1}(P, \reals )$,
which will be  called the {\em flux homomorphism}.
To define   $\frakf$ precisely, we will  first embed $U(\gm )$ into
  the symplectomorphism group $Diff(\gm ,  \omega )$ by  considering each
lagrangian bisection as  a symplectic diffeomorphism on
$\gm$ via left translations, and then take the usual flux homomorphism
for  symplectomorphisms.
An  advantage
of this approach is that 
one can directly apply  various   results about  the usual flux homomorphism
and avoid  many  tedious  verifications.
Discussion on  the flux homomorphism and  its
properties  occupies Section 2 and Section 3.
In the case that $P$ is symplectic, 
 the symplectic structure establishes
an isomorphism between  the cotangent bundle  $T^{*}P$
and the tangent bundle $TP$.
Hence the space of one-forms can be identified
  with the space of vector fields. Then $Z^{1}(P)$ and $B^{1}(P)$
are isomorphic, as  Lie algebras, to the space of
symplectic vector fields and  that  of hamiltonian vector
fields,  respectively. Moreover, the pair groupoid 
$\gm =P\times \Bar{P}$ can be taken as a symplectic
groupoid over $P$. The group $U(\gm )$   is thus
isomorphic to the symplectomorphism
group $Diff (P, \omega )$, and $\frakf$ reduces to the usual flux
homomorphism.

To lift the  homomorphism $ C^{\infty}(P)\stackrel{d}{\lon } Z^{1}(P)$,
first of all, we need a Lie group integrating the Lie
algebra $C^{\infty}(P)$.  It turns out that such  a Lie group can be described
 in terms of  a prequantum bundle $E$ over the symplectic groupoid  $\gm$ 
(properties of such  prequantum bundles 
 were  studied in \cite{WeinsteinX:1991}).
We will devote Section 4 to  the discussion on  this aspect.

Despite its conceptual clearness, the flux homomorphism
is in general very difficult to compute.
In Section 5,  
we will study  an  alternate map
closely related to  the flux homomorphism   when there exists
 an invariant measure $\nu$ on the  Poisson manifold $P$.
This is   the composition $\mu \smalcirc \frakf :
\widetilde{\gmm}\lon H^{1}(P, \reals )^{*}$ of the
flux homomorphism $\frakf$ with a  linear map
$\mu :  H^{1}(P, \reals )\lon  H^{1}(P, \reals )^{*}$ 
induced from the  invariant measure $\nu$.
It turns out that this  map  $\mu \smalcirc \frakf$
can be   explicitly  expressed  by using integration
when  $\gm$ is $\alpha$-simply
connected. This   becomes especially  useful
when $\mu$ is an isomorphism. In this case,
one    can  indeed     compute the flux 
homomorphism by  computing  $\mu \smalcirc \frakf$.
As an example, this method will be applied 
to symplectic manifolds in Section 6.

For simplicity,  we  will restrict our attention  to
compact    Poisson manifolds.
  However, many   results should  be easily  modified for noncompact
 Poisson manifolds when 
 various notions in the context are replaced 
 by their compactly supported analogues.

Recently, it came to the author's attention
 that the  flux homomorphism (in a more general  setting)
  was also been constructed
 by Dazord   \cite{Dazord:1994}.
However, our approach  is  quite different.
While Dazord   deals with  general noncompact Poisson manifolds with 
possibly non-Hausdorff symplectic groupoids, he 
 proves   his assertions directly.
We  will instead  work within the framework of compact Poisson manifolds 
admitting  Hausdorff symplectic groupoids.
  Within  this  more conventional framework, it  allows us to apply  
many  known results  about  the usual
flux homomorphism of symplectomorphism groups \cite{Banyaga}
\cite{McduffS}. Since the original preprint, another
related paper \cite{Dazord:1995} has also appeared.

{\bf Acknowledgements.} The author  would like to thank
 Dmitry Fuchs, Dusa McDuff, Lisa Traynor  and Alan Weinstein 
for useful discussions.
He is particularly grateful to Dusa McDuff
for letting  him have access  to the  unpublished
manuscript of a book by  her with   Salamon.
Part of this work  was finished when the author was
visiting Isaac Newton Institute in the summer of 1994.
He wishes to thank the Institute for its hospitality and
financial support. Special thanks goes to Augustin  Banyaga 
for bringing his attention to the paper \cite{Dazord:1994}.
 Finally, he also wishes to thank the referees
for their helpful comments and suggestions.
 
\section{The Flux homomorphism}

Let  $\poiddd{\gm}{P}{}$ be a   symplectic groupoid over a
compact Poisson manifold, with 
symplectic structure $\omega$.  A lagrangian
submanifold of $\gm$ for which the restrictions of $\alpha$ and
$\beta$ are diffeomorphisms onto $P$ is called a lagrangian bisection.
Under the multiplication of subsets induced from the product on $\gm $,
the lagrangian bisections form a Lie group $U(\gm)$ in terms
of Whitney $C^{\infty }$-topology.
$U(\gm)$ is easily seen to be locally path connected by the
local structure theorem of lagrangian submanifolds \cite{Weinstein:1971}.
The unit of  $U(\gm)$ is the identity
space $P$ of the groupoid,  and is   denoted by $1$.
 We  denote
by $\gmm $ the identity component of $U(\gm)$. For any $L\in U(\gm)$
and $u\in P$, $L(u)$ denotes the unique point in $L$ characterized by
 the condition
that $\beta L(u)=u$. A lagrangian bisection $L$
induces a Poisson diffeomorphism  on  the base Poisson manifold $P$  by
$u\lon \alpha (L(u))$, which will be  denoted by $Ad_{L}(u)$.

Let $L_{t}$ be an isotopy  in $\gmm$ starting at  the identity
section, i.e., $L_{0}=1$.  For any fixed $t$, the variation of $L_{t}$
around the neighborhood of that fixed time   defines a  vector
field $X_{t}$ along the bisection $L_{t}$ on $\gm $.
That is,

\begin{equation}
\label{eq:X}
X_{t}(L_{t}(u))=\diff{s}|_{s=t}L_{s}(u).
\end{equation}

It follows from the property $\beta L_{t}(u)=u$ that
$$\beta_{*}X_{t}=0. $$
That is, $X_{t}$ is tangent to $\beta $-fibers. This implies that
there exists a time-dependent  one-form $\theta_{t}\in \Omega^{1}(P)$ 
on the base such that

\begin{equation}
\label{eq:one-form}
X_{t}\per \omega =\alpha^{*} \theta_{t} 
\end{equation}
along the bisection $L_{t}$.

 For any $L\in  U(\gm )$, we denote by  $\tau (L)$
 the left translation  on $\gm $ by $L$, which is a diffeomorphism 
of $\gm$.  In fact, $\tau (L)$ is  a  symplectomorphism of  $\gm$.
We thus obtain a map
 $\tau : U(\gm )\lon Diff(\gm , \omega )$, which is a group
homomorphism.  Clearly, $\tau$ maps
the identity component $ \gmm $ into $Diff_{0}(\gm , \omega )$. 
By $\tilde{\tau}:\widetilde{ U(\gm)_{0} }
\lon \widetilde{ Diff_{0}(\gm , \omega ) }$, 
we denote its  natural lifting  to  
their universal coverings.
 
 Let $\phi_{t} =\tau (L_{t}) $ be the  family of
symplectomorphisms of  $\gm$
corresponding to  the   lagrangian bisection  isotopy  $L_{t}$.
We denote by $Y_{t}$   the 
time-dependent vector field on $\gm $
which generates  $\phi_{t}$, i.e., 

\begin{equation}
\label{eq:Y}
\diff{t} \phi_{t}=Y_{t}\smalcirc \phi_{t}.
\end{equation}

The usual flux homomorphism of  $\{\phi_{t} \}$, considered as
a homotopy class of a   symplectic  isotopy of  $\gm$,
is defined in terms of  such  a time-dependent  vector field 
$Y_{t}$ \cite{Banyaga} \cite{Calabi}. To find out  the precise relation
between  $Y_{t}$ and  $\theta_{t}\in \Omega^{1}(P)$,
we first need to compute the one-form $ Y_{t}\per \omega $.
\begin{pro}
\label{pro:Y}
The time-dependent one-form $Y_{t}\per \omega $ on $\gm $ is equal to
$\alpha^{*} \theta_{t} $.
\end{pro}
\pf  Fix any $t=\ttt$,

\be
Y_{\ttt}(\phi_{\ttt}(x))=&=&\diff{t}|_{t=\ttt}\phi_{t}(x)\\
&=&\diff{t}|_{t=\ttt}(L_{t}(\alpha (x)) \cdot x)\\
&=&(r_{x})_{*} X_{\ttt}|_{L_{\ttt}(\alpha (x))},
\ee
where the right hand side makes sense since $X_{\ttt}$ is
tangent to $\beta$-fibers. 

Let $\calx$ be any lagrangian bisection through the point $x$. Then,
\be
Y_{\ttt}(\phi_{\ttt}(x))\per \omega &=& [ (r_{\calx })_{*} X_{\ttt}|_{L_{\ttt}
(\alpha (x))} ]\per \omega\\
&=&(r_{\calx })_{*} [X_{\ttt}|_{L_{\ttt} (\alpha (x))}\per \omega ]\\
&=&(r_{\calx })_{*} \alpha^{*}[\theta_{\ttt}|_{Ad_{ L_{\ttt}}\alpha (x)}]\\
&=&\alpha^{*}[\theta_{\ttt}|_{Ad_{ L_{\ttt}}\alpha (x)}].
\ee

By letting $y=\phi_{\ttt} (x)$,
one gets immediately that

$$Y_{\ttt}(y)\per \omega =\alpha^{*}[\theta_{\ttt}(\alpha (y))]. $$
This concludes our proof.
\qed

An immediate consequence is the following
\begin{cor}

$\theta_{t}$ is closed for all $t$.
\end{cor}
\pf It follows from $L_{Y_{t}}\omega =0$ that $Y_{t} \per \omega $
is closed for all $t$.  This  implies immediately
that $\theta_{t}$ is closed. \qed

As a by-product, below we will see   how   the Lie algebra of $U(\gm )$
can be  identified with the Lie algebra $Z^{1}(P)$
of closed one-forms   for compact Poisson manifolds
\cite{Weinstein:1987}.
Since we need to use such an identification frequently
in the paper, we briefly
recall it below.
Given a smooth path of lagrangian bisections $\{ L_{t }\}$  starting
from $1$,
its  derivative at   $t=0$
  is given by the section 
$u\lon \diff{t}|_{t=0}L_{t}(u)$ of the vector bundle  $T_{P}^{\beta}\gm$,
which is isomorphic to the normal bundle $T_{P}\gm /TP$.
Being  contracted  with the symplectic form  $\omega $, 
this section  can be naturally identified with a one-form
on $P$.
This  is exactly  the closed one-form $\theta_{0}$ as introduced
earlier in this section. In this way, we obtain an identification
of the tangent space of $U (\gm )$ at $1$ with  the
space $Z^{1}(P)$ of closed one-forms.

We now turn to the discussion on the flux homomorphism.
It follows from Proposition \ref{pro:Y} that the flux homomorphism of
$ \{ \phi_{t} \} $ is given by

\begin{equation}
\label{eq:flux-rel}
F( \{\phi_{t} \} )=\alpha^{*} \int_{0}^{1} [\theta_{t}]dt .
\end{equation}
This suggests the following 

\begin{defi}
We define the flux homomorphism on a symplectic groupoid
$\gm $ to be the map $\calf :\widetilde{\gmm }\lon H^{1} (P, \reals )$
given by:

\begin{equation}
\frakf  (\{ L_{t} \} )=\int_{0}^{1} [\theta_{t}]dt, \ \ 
\mbox{for  any \ }  \{ L_{t} \}\in \widetilde{\gmm }.
\end{equation}
\end{defi}

\begin{pro}
Suppose that $P$ is  a compact Poisson manifold.
Then, $\frakf$ is  well-defined  and  surjective.
Moreover,   it is a  group homomorphism when $H^{1}(P, \reals )$
 being  considered as an Abelian
group.
The  derivative of $\frakf$ is the map 
$f: \ Z^{1}(P) \lon H^{1}(P, \reals )$ which sends 
 any closed 1-form to
its cohomology class. 
\end{pro}
\pf It follows 
from   Equation (\ref{eq:flux-rel}) that the following diagram:

\begin{equation}                         \label{tdvb}
\matrix{&& \tilde{\tau} &&\cr
        &  \widetilde{ U_{0}(\gm)} &\vlra&  \widetilde{ Diff_{0}(\gm , \omega ) }&\cr
        &&&&\cr
     \frakf &\Bigg\downarrow&&\Bigg\downarrow&F \cr
        &&&&\cr
        & H^{1}(P, \reals ) &\vlra&  H^{1}(\gm , \reals )&\cr
        && \alpha^* &&\cr}
\end{equation}
commutes.

Let $i:P\lon \gm$ be the inclusion map. It follows from
$\alpha \smalcirc i=id$ that $i^{*}\smalcirc \alpha^{*}=id$, which 
 implies that $\alpha^* :  H^{1}(P, \reals ) \lon H^{1}(\gm , \reals ) $
is injective.
Since $F$ is  well-defined on homotopy classes of symplectomorphisms,
it follows that $\frakf$ is  also well-defined.  The
same argument shows that $\frakf$
is  indeed a homomorphism.

 To prove the second assertion, 
 we consider $L_{t}=\exp{t\theta}$, $0\leq t\leq 1$, 
 for any closed  one-form $\theta $.
Here by definition, $\exp{t\theta}$ is a family of
lagrangian bisections obtained  by moving the identity  space $P$
along the flows of $X_{\alpha^{*}\theta }$ \cite{M:B1} \cite{Xu:OPG}.
Clearly, $L_{t}$ is defined for all $t$ since $X_{\alpha^{*}\theta }$
is a complete vector field.
It is simple to see 
 that $\frakf \{ L_{t} \} =[\theta ]$. The same argument also
shows that the derivative of $\frakf$ is indeed given
by  the map: $\theta \in Z^{1}(P  )\lon [\theta] \in  H^{1}(P, \reals )$.
This concludes the proof of the proposition. \qed

We end this section with some simple examples.

\begin{ex}
If $P$ is a compact symplectic manifold, the pair groupoid
$\gm =P\times \Bar{P}$ can be taken as  its symplectic groupoid.
Any  lagrangian bisection of $\gm$
can be naturally 
 identified with the graph of a symplectic diffeomorphism.  Therefore,
$U(\gm )$ can be identified with  the group of symplectomorphisms.
The flux homomorphism, in this case,
 reduces to the usual flux homomorphism 
of symplectomorphisms.
\end{ex}

\begin{ex}
If $P$ is a compact manifold with zero Poisson structure,
 the symplectic groupoid  $\gm $ is the cotangent space $T^{*}P$.
 Then $U(\gm )$ is the space of all closed one-forms, which is
connected and simply connected. The group multiplication is
simply the addition of one-forms.
The flux homomorphism $\frakf : U_{0}(\gm )\lon H^{1}(P, \reals )$
is then  the map sending any closed one-form to its cohomology class.
\end{ex}

\begin{ex}
Let $P$ be  a regular Poisson manifold which is the product
of a zero Poisson manifold $Q$ with a symplectic
manifold $S$. Its  symplectic groupoid is then 
the product groupoid $T^{*}Q\times S\times \Bar{S}$. 
A bisection of $\gm $ consists of a pair $\phi \in C^{\infty}(S, \Omega^{1}(Q))$
and $\psi \in C^{\infty}(Q, Diff(S))$. If it  is a lagrangian
bisection, its restrictions  to both  the inverse images
  $\alpha^{-1}(\{q\}\times S)$ 
and $\alpha^{-1}(Q\times \{x\})$ are isotropic submanifolds for any
$q\in Q$ and $x\in S$.
This implies that the image of $\psi $ should lie in the symplectic
diffeomorphism group $Diff (S, \omega )$ and  $\forall x\in S$,
$$d\phi (x)+\psi_{x}^{*}\omega =0, $$
where $\psi_{x} : Q\lon S $ is the evaluation map at $x$
defined by $\psi_{x}(q)=\psi (q) \cdot x$.
Therefore, $U(\gm )$ is a subgroup of
$C^{\infty}(S, \Omega^{1}(M))\times C^{\infty}(M, Diff(S, \omega ))$. 
The explicit expression of the flux
homomorphism is rather involved and we shall omit it here.
\end{ex}

\section{Exact lagrangian bisections}

The Lie algebra of all closed one-forms contains the subspace
$B^{1}(P)$
of all exact one-forms  as a  Lie subalgebra. It is therefore
expected that  $\gmm $  contains  a subgroup
having    $B^{1}(P)$ as its Lie algebra.  The purpose of the
present section is to characterize  the  elements in
this  subgroup by using  the flux homomorphism introduced in the  previous
section.

\begin{defi}
 A lagrangian bisection $L$  in $ \gmm $ is said to be exact iff
there is  an isotopy $L_{t}\in U(\gm )$ connecting $1$ and $L$ (i.e., $L_{0}=1$ 
and $L_{1}=L$) such that $\tau (L_{t})$ is a hamiltonian isotopy
of $\gm$.
\end{defi}

When $P$ is symplectic and $\gm$  is the pair  groupoid
$P\times \Bar{P} $, a lagrangian bisection of $\gm$ 
corresponds to a symplectomorphism on $P$.
It is simple to see that  exact  lagrangian bisections
  will correspond  to hamiltonian symplectomorphisms.

We denote by $\gmex$ the space of all exact lagrangian bisections.
Since the product of two  hamiltonian isotopies is still a
hamiltonian isotopy, $\gmex$ is a subgroup of $ \gmm $.

\begin{pro}

$\gmex$ is a path-connected normal subgroup of $\gmm$.
\end{pro}
\pf This follows from the fact that the group of   hamiltonian  symplectomorphisms 
$Ham (\gm , \omega )$  
is a normal subgroup of  $Diff_{0}(\gm , \omega )$.
\qed

The following proposition is an immediate consequence of
Proposition \ref{pro:Y}.

\begin{pro}

A lagrangian bisection $L\in \gmm$ is exact iff there is
an isotopy  $L_{t}\in U(\gm )$ connecting $1$ and $L$ 
such that the corresponding time-dependent one-form $\theta_{t}$ defined by
Equation (\ref{eq:one-form})   is   exact for all $t$.  
\end{pro}

The same argument also shows the following:

\begin{pro}
For a compact Poisson manifold $P$, 
the Lie algebra of $\gmex$ is  the space $B^{1}(P)$
of all exact one-forms on $P$.
\end{pro}
\pf Let $L_{t}$ be any path in $\gmex$ starting from $1$, with
$\theta_{t}$ being the corresponding time-dependent one-form  on $P$.
By definition,  $\tau (L_{t})$ is a hamiltonian isotopy in $\gm$,
so $\alpha^{*} \theta_{t}$ is exact for all $t $ according
to Proposition \ref{pro:Y}. Therefore,  $ \theta_{t}$ is 
all exact, so in particularly is $\theta_{0}$.
Conversely, given any exact one-form $\theta$ on $P$,
$\exp{t\theta }$ will be  a path in $U_{ex}(\gm  )$.
This shows that the tangent space of $U_{ex }(\gm )$
at $1$ is isomorphic to $B^{1}(P)$.
\qed

\begin{pro}
\label{pro:path-flux}
Assume that the base Poisson manifold $P$ is compact.
Let $L_{t}$ be a path in $U ( \gm )$ such that $L_{0}=1$ and
$L_{1}=L$. Then,  $L_{t}$ lies in $\gmex$ for
all $t \in [0, 1]$ if and only if $\frakf (\{ L_{t} \}_{0\leq t\leq T})=0$,
for all $0\leq T \leq 1$.
\end{pro}
\pf By definition,  that $L_{t}$ is  in $\gmex$,  for all $t$,  is
equivalent to saying that $\tau (L_{t})$  lies in 
 $Ham (\gm , \omega )$, for all $t$.  The latter is
equivalent to that $\tau (L_{t})$  itself is  a hamiltonian
isotopy, which is equivalent to $F (\{\tau  L_{t} \}_{0\leq t\leq T})=0$,
$\forall 0\leq T\leq 1 $, 
according to  Proposition II 3.3 in \cite{Banyaga}
 or Theorem 5.2.4 in \cite{McduffS}. Now $F (\{\tau  L_{t} \}_{0\leq t\leq T})
 = \alpha^{*}\frakf (\{ L_{t} \}_{0\leq t\leq T})$. The conclusion thus
follows from the fact that $\alpha^*$ is injective. \qed
{\bf Remark\ }  Either   Proposition II 3.3 in \cite{Banyaga}
or  Theorem 5.2.4  in \cite{McduffS}  requires that
the symplectic manifold be compact.
In our case, the symplectic groupoid 
$\gm$ is generally
not compact even though $P$ is assumed to be compact.
However,  we note that the requirement for the compactness  of
the symplectic manifolds  arises from the 
 necessity   that certain relevant   vector fields  be complete.
In our  case,
the completeness of these  vector fields is guaranteed
by the compactness of $P$ according to
 a result  of Kumpera and Spencer (see
the lemma in  Section 33 of Appendix of  
 \cite{KumperaS}).

By using the flux homomorphism, we now can characterize 
exact lagrangian bisections.

\begin{thm}
Assume  that the base Poisson manifold  $P$  is compact.
Then a  lagrangian bisection $L$ is exact 
 in $\gmex$ if and only if
there is an isotopy $L_{t}$ in $U(\gm )$ connecting $1$ and
$L$ such that $\frakf (\{L_{t}\})=0$.
\end{thm}
\pf Our proof here is essentially borrowed  from that of Theorem 5.2.4 
 in \cite{McduffS}.
One direction is quite straightforward by  definition.
Namely, if $L$ is exact, then there  exists an isotopy
 $L_{t}$  in $\gmex$ connecting  $1$ and $L$.
Thus, it follows from Proposition \ref{pro:path-flux}
 that $\frakf (\{L_{t}\})=0$.

Conversely, assume that $L_{t}$ is an isotopy  in $U(\gm )$ connecting $1$ and
$L$  such that $\frakf (\{L_{t}\})=0$. Let
$\theta_{t}$ be the corresponding time-dependent closed   one-form on $P$
as defined by Equation (\ref{eq:one-form}).
Then,  $\int_{0}^{1} \theta_{t} dt$ is exact.

First, we assume that $\int_{0}^{1} \theta_{t} dt=0$.

Let $\gamma_{t}$ be a family of  closed one-forms  on $P$ defined  by
$\gamma_{t}=-\int_{0}^{t} \theta_{s} ds$, $0\leq t\leq 1$ and
$K_{t}=L_{t}\exp{\gamma_{t}}, \  0\leq t\leq 1$. 
Here, the map $\exp{}$ is well defined since $P$ is
compact.
Then
$K_{t}$ is an isotopy of   lagrangian bisections connecting
$1$ and $L$. Moreover, for any $T\in [0, \ 1]$, 

\be
\frakf (\{K_{t}\}|_{0\leq t\leq T} )&=&
\frakf (\{L_{t}\}|_{0\leq t\leq T} )+\frakf  ( \{\exp{\gamma_{t}}\}|_{0\leq t\leq T} )\\
&=&[\int_{0}^{T} \theta_{t} dt] +\frakf (\{\exp{t\gamma_{T}}\}|_{0\leq t\leq 1
} )\\
&=&[\int_{0}^{T} \theta_{t} dt] +[\gamma_{T}]\\
&=&0,
\ee
where  the second equality follows from the homotopy
equivalence between 
 the two isotopies  $\{ \exp{\gamma_{t}}\}|_{0\leq t\leq T} $ and
$\{\exp{t\gamma_{T}}\}|_{0\leq t\leq 1}$. This shows that $K_{t}$ is
a path in $\gmex$. In particular, $K_{1}$ is exact, which implies
that $L$ must be exact.

Next, let us  assume that $\int_{0}^{1} \theta_{t} dt=d H$ for 
some function $H\in C^{\infty}(P)$.
Consider $L_{t}'=\exp({-tdH })*L_{t}$, where the $*$-product 
means the following:

$$ \exp({-tdH })*L_{t}=\left\{ \begin{array}{ll}
                   L_{2t} & 0\leq t\leq \half\\
                   \exp({-(2t-1) dH)} \cdot L_{1}, & \half \leq t \leq 1.
\end{array}
\right.  $$

Let  $\theta_{t}'$  denote 
its corresponding time-dependent  one-form on $P$
as defined by Equation (\ref{eq:one-form}). Then,
it is clear that 

$$
\theta_{t}'=\left\{ \begin{array}{ll}
                   2 \theta_{2t} & 0\leq t\leq \half\\
                  -2 dH , &  \half  \leq t \leq 1 .
\end{array}
\right.
$$
Therefore,  $\int_{0}^{1} \theta_{t}' dt=0$.
Hence, $\exp({-dH})\cdot  L\in \gmex$,  and then  $L\in \gmex$. \qed

As an immediate consequence, we have

\begin{cor}
\label{cor:ex-cover}
Suppose that $P$ is compact.
The kernel of the flux homomorphism $\frakf :\widetilde{\gmm}\lon H^{1}(P, \reals )$
is the universal covering of $\gmex$. I.e., $ker \frakf = \widetilde{\gmex}$.
\end{cor}
\pf  It is clear that $\widetilde{\gmex}\subseteq ker \frakf$.
To prove the other half, assume  that $L_{t}$ is an isotopy of
lagrangian bisections   satisfying 
$\frakf ( \{L_{t }\})=0$. Let $\theta_{t}\in \Omega^{1}(P)$ be
its   corresponding time-dependent  one-form on $ P$.
Then, we have   $\int_{0}^{1} \theta_{t} dt=d H$ for some function $H$ on $P$.
According to the  proof  of the  last proposition, 
$\exp{(-tdH) }*L_{t}$ is  an isotopy in $\gmex$.
On the other hand, it is easy to see that $L_{t}$ is
homotopy equivalent to $\exp{(tdH)}*[\exp{(-tdH )}*L_{t}]$.
The latter is clearly a path entirely in $\gmex$, so it follows that
$\{ L_{t} \}\in  \widetilde{\gmex}$. \qed
{\bf Remark}  For symplectic manifolds, the flux homomorphism is useful
because it can be used to characterize hamiltonian symplectomorphisms.
Hamiltonian symplectomorphisms,  rather than symplectomorphisms
are what Arnold's conjecture concerns about. It would be interesting to ask
if there is something similar to
 Arnold's conjecture for Poisson manifolds, which should 
concern about the  intersection  number of
 an exact lagrangian bisection with
the identity space. If the Poisson manifold is a compact zero Poisson
manifold, this would fall into the content of Morse theory.
However, for a  general Poisson manifold,
it is not even clear whether an exact lagrangian 
bisection always intersects with the identity space,
although there are some evidences that the answer 
 is likely positive.

\section{A Central extension  of $\gmex$ }

In this section, we will assume that the symplectic groupoid 
$(\gm , \omega ) $ is prequantizable. That is, $\omega$ is  of 
integer class.  
In this case,  the group $\gmex$ can be
easily characterized with  the help of the prequantum bundle over $\gm$.
More importantly, such a prequantum bundle provides us with
a natural  one-dimensional central extension of $\gmex$, whose corresponding
Lie algebra is    the Poisson algebra $C^{\infty}(P)$.
Thus,   we will  obtain   a   geometric description for a  Lie group
integrating the   Poisson  
 algebra $C^{\infty}(P)$.

Suppose that $\gm$ is prequantizable  so that we can construct
a unique circle ($\reals/\integers$)-bundle $p: E\lon \gm$ with a  
connection form $\theta $ satisfying the condition that the
identity space $P$ has no holonomy. It is proved in \cite{WeinsteinX:1991}
that $E$ has a groupoid structure once a horizontal lifting of
$P$  is chosen in $E$. In the sequel, we will fix   such  a horizontal
 lift.
 Thus, $E$ carries  a groupoid structure.  As in \cite{WeinsteinX:1991},
we will denote the source   and target maps of $\poidd{E}{P}$  by
$\alpha$ and $\beta$, while those of $\gm $ are
denoted by $\alpha_{0}$ and $\beta_{0}$.
Then, $\alpha=\alpha_{0}\smalcirc p$ and $\beta=\beta_{0}\smalcirc p$.

For a lagrangian bisection $L\in U(\gm )$, the restriction of $E$
on  $L$ is  always flat.
For  any loop $\gamma (t)$ in $P$, $L \cdot \gamma (t)$
is a loop in $L$.
We define $\Phi_{L}(\gamma )$ as its corresponding holonomy.
Thus, we obtain a map $\Phi :L\lon \Phi_{L}$
from $ U(\gm )$ to $Hom (\pi_{1}(P) , T^{1})$. Here 
$\pi_{1}(P)$ denotes the fundamental group of $P$.
It is clear that  $Hom (\pi_{1}(P) , T^{1})$  is a left  $U(\gm )$-module
with the action:  $(L\cdot \phi )(\gamma )=\phi (Ad_{L} \gamma )$,
where $L\in U(\gm )$,    $\phi \in Hom (\pi_{1}(P) , T^{1})$ and $[\gamma ]
\in \pi_{1}(P )$.

\begin{pro}

The map $\Phi: U(\gm )\lon Hom (\pi_{1}(P) , T^{1})$  satisfies
the following 1-cocycle-like property:
$$\Phi_{KL}=\Phi_{L}+L\cdot \Phi_{K}. $$
\end{pro}
\pf Suppose that  $\gamma (t)$ is  any loop in $P$. Let $ \gamma' (t)=
Ad_{L}   \gamma (t)$, $\gamma_{1}(t)=L  \cdot \gamma (t)$ and 
$\gamma_{2}(t)=K  \cdot \gamma' (t)$. Then it is clear that 
$\beta \gamma_{2}(t)=\alpha \gamma_{1}(t)= \gamma' (t)$ and
$\gamma_{2}(t)  \gamma_{1}(t)=KL\cdot \gamma (t)$.
Therefore,  $( \gamma_{2}(t),  \gamma_{1}(t), \gamma_{2}(t)  \gamma_{1}(t))$ is
a loop lying in the graph of groupoid multiplication
 $\Lambda \subset \gm \times \gm \times \Bar{\gm}$.
According to  Theorem 3.1 in \cite{WeinsteinX:1991}, 
$\Lambda$ has no holonomy in the corresponding prequantum
bundle $(E\times E\times \Bar{E})/T^{2}\lon \gm \times \gm \times \Bar{\gm}$.
It thus follows that $Hol( \gamma_{2}\gamma_{1})=Hol (\gamma_{2})+
Hol(\gamma_{1})$. I.e., $\Phi_{KL}(\gamma )= \Phi_{K}(Ad_{L}\gamma )+\Phi_{L}(\gamma )$.
\qed

An immediate consequence  is the following:

\begin{cor}

$\Phi: \gmm\lon Hom (\pi_{1}(P) , T^{1})$ is a group homomorphism.
\end{cor}

We will see below that 
 $\Phi$ is closely
related to the flux homomorphism. To describe their  precise relation, 
let us introduce a  group homomorphism $\rho: 
H^{1}(P, \reals )\lon Hom (\pi_{1}(P) , T^{1})$ 
by 
$$\rho ([\theta]) ([\gamma ] )=\int_{\gamma }\theta \ \ \ \mbox{ mod}
 \ \integers, \  \ 
\forall [\theta ]\in H^{1}(P , \reals ),  \mbox{  and }
[\gamma ]\in \pi_{1} (P) .$$

\begin{pro}
\label{pro:flux-phi}
Suppose that $P$ is compact.
For any isotopy  $L_{t}\in \gmm$ connecting $1$ and $L$,
$$(\rho \smalcirc \frakf )( \{L_{t}\} )=\Phi (L). $$
\end{pro}
\pf Suppose that $\gamma (s),\ 0\leq s\leq 1$ is any 1-cycle  in $P$.
Let $\lambda (t, s)=L_{t}^{-1} \gamma (s)$, and
$\theta_{t}$ be the time-dependent   closed 1-form on $P$  corresponding to
 $L_{t}$. 
As before, $Y_{t}$ denotes  the time-dependent vector field on $\gm$
corresponding to the  one family of diffeomorphisms  $\phi_{t}=\tau (L_{t})$
defined by Equation (\ref{eq:Y}). 

 Taking the   derivative   of   the identity $\phi_{t}(L_{t}^{-1} x)=x$
 at $t=t_{0}$, one gets:
$$Y_{t_{0}}(x)+(\phi_{t_{0}})_{*} \diff{t}|_{t=t_{0}}L_{t}^{-1} x=0.$$
Since $x$ is an arbitrary point in $\gm$, by taking $x=\gamma (s) L_{t_{0}}$,
 we have 
$$Y_{t_{0}}( \gamma (s) L_{t_{0}} )
+(\phi_{t_{0}})_{*} \diff{t}|_{t=t_{0}}L_{t}^{-1} \gamma (s) L_{t_{0}} =0.$$
I.e., 
\begin{equation}
Y_{t_{0}}( \gamma (s) L_{t_{0}} ) =-L_{t_{0}} \pttt{\lambda (t, s)}{t}|_{t=t_{0}} L_{t_{0}}.
\end{equation}
Therefore,
\be
\dot{\gamma}(s)\per \theta_{t}&=& \dot{\gamma}(s)L_{t}\per \alpha^{*}  \theta_{t}\\
&=&
\omega (Y_{t}(\gamma (s) L_{t}),  \dot{\gamma}(s)L_{t})\\
&=&\omega (-L_{t} \pttt{\lambda}{t} L_{t}, L_{t} \pttt{\lambda}{s} L_{t})\\
&=&-\omega (\pttt{\lambda}{t} ,  \pttt{\lambda}{s} ), 
\ee 
where  the last step follows from the fact that $L_{t}$ is
lagrangian.
Hence,
\be
<\rho (\frakf \{L_{t}\}), [\gamma]>&=&\int_{\gamma}\int_{0}^{1}\theta_{t}dt\\
&=& \int_{0}^{1}\int_{0}^{1} \dot{\gamma}(s)\per \theta_{t} \  dt ds\\
&=&-\int_{\cal D} \omega,
\ee
where ${\cal D} $ is the  cylinder: $\lambda (t,s)=L_{t}^{-1}\gamma (s)$,
$0\leq t\leq 1$ and $s\in S^{1}$.
Now the quantity: $\int_{\cal D} \omega$ modulo  $\integers$
 is exactly  the difference between the   holonomies around
 the loops  $\lambda (1, s)$ and
$\lambda (0, s)$. The holonomy of the latter is zero because the loop
lies in $P$, while the  holonomy of $\lambda (1, s)$ is the minus
of that of $L \cdot \gamma (s)$, which is  $\Phi_{L}(\gamma (s))$.
This completes the proof of the proposition. \qed

We denote by $\Omega$ the image of the fundamental group
$\pi_{1}(\gmm , 1)$ under the flux homomorphism $\frakf$, i.e.,
$\Omega=\frakf (\pi_{1}(\gmm , 1))$. It follows from  
 Proposition \ref{pro:flux-phi}
 that $\rho (\Omega )=0$. Therefore, $\rho$ descends
to a homomorphism:  $H^{1}(P, \reals )/\Omega \lon  Hom (\pi_{1}(P), T^{1})$,
which will  still be  denoted by $\rho$.
We note that the flux homomorphism induces a   homomorphism:
$\gmm\lon H^{1}(P, \reals )/\Omega$, which will be denoted as $\tilde{\frakf}$
below.

The following theorem is an alternate description of 
 Proposition \ref{pro:flux-phi}.

\begin{thm}
\label{thm:flux-phi}
Assume that $P$ is a compact  Poisson manifold and its
 symplectic groupoid $\gm $ is prequantizable.
Then, the following diagram:
\begin{equation}                         \label{flux-hol}
\matrix{&& \tilde{\frakf } &&\cr
        &   \gmm &\vlra&   H^{1}(P, \reals )/\Omega
 &\cr
        &&&&\cr
     \Phi  &\Bigg\downarrow&&\Bigg\downarrow&\rho  \cr
        &&&&\cr
        &  Hom (\pi_{1}(P), T^{1}) &\vlra&   Hom (\pi_{1}(P), T^{1} )&\cr
        && id &&\cr}
\end{equation}
commutes.

\end{thm}

Using this result, we can  describe the group $\gmex$ in
terms of  $\Phi$.

\begin{cor}
\label{cor:component}
Under the same hypotheses of Theorem \ref{thm:flux-phi},
 the group $\gmex$ coincides with the identity component of
the kernel of $\Phi$. I.e.,
$$\gmex=(Ker \Phi )_{0}. $$
\end{cor}
\pf For any $L$ in $\gmex$, by definition, there exists
an isotopy  $L_{t}$ connecting $1$ and $L$ such that 
$\tau (L_{t} )$ is a hamiltonian isotopy.
Hence, $\frakf (\{L_{s}\}_{0\leq s\leq t})=0$ for any
$t\in [0, 1]$. In other words, $\tilde{\frakf}(L_{t})=0$ for all $t\in [0, 1]$.
It follows from Theorem \ref{thm:flux-phi} that $\Phi (L_{t})=
(\rho  \tilde{\frakf}) (L_{t})=0$. Hence, $L_{t}$ is a path in $Ker \Phi$.
This   implies that
$L\in(Ker \Phi )_{0}$.

Conversely, suppose that $L_{t}$,  with $L_{0}=P$ and $L_{1}=L$,  is
an isotopy  in $ Ker \Phi $. Thus,
$(\rho \tilde{\frakf})(L_{t} )=\Phi (L_{t})=0$. This implies that
for any loop $\gamma (s)$ in $P$, the double integral
$\int_{\gamma }\int^{t}_{0} \theta_{s}ds $ is an integer for all $t$. 
However, since  it is a continuous function with respect to $t$,
it must  be a constant and in our case
it should be  identically zero.
Therefore,  $\frakf( \{L_{s}\}_{0\leq s\leq t} )=0, \ \forall t$.
Hence, $L$ is in $\gmex$. \qed

An immediate consequence is the following:

\begin{cor}
Under the same hypotheses, the Lie algebra of $Ker \Phi $ is $Z^{1} (P)$.
\end{cor}

The prequantum bundle $E$ together with the connection
form $\theta$ becomes a contact manifold. A bisection
of the groupoid $\poidd{E}{P}$ is called a legendre
bisection if it is also a legendre submanifold
with respect to this  contact structure (i.e. horizontal
with respect to the connection $\theta $).
It is clear that the space of all  legendre  bisections of
 $\poidd{E}{P}$ is a Lie  group.
We denote by $U_{0}(E)$  its connected component of
the unit.
It is a subgroup of the group of all bisections of   $\poidd{E}{P}$.

It is well known that  left (or right)
 translations by lagrangian bisections are symplectic diffeomorphisms of 
the symplectic groupoid.
Similarly,  left (or right) translations by legendre  bisections
 are
 contact diffeomorphisms.

\begin{pro}
\label{pro:contact}
A bisection $K$ of   $E$ is legendre  iff  the 
left translation $l_{K}$ (or the right translation $r_{K}$): $E\lon E$
is a contact diffeomorphism.
\end{pro}
\pf Let $\psi$ denote the diffeomorphism $l_{K}: E\lon E$. 
For any tangent vector $\delta_{x}$ in $E$, we take a
path $x(t)$  starting from $x$ with $\dot{x} (0) =\delta_{x}$.
Let $y(t)=K(\alpha (x(t)))$. Then $y(t)$ is a path in $K$
and $(y(t), x(t), \psi (x(t))) $ lies in the groupoid  graph of $E$.
Since the one-form
 $(\theta , \theta , -\theta)$ vanishes on the groupoid  graph of $E$ \cite{WeinsteinX:1991}, 
and $K$ is a legendre submanifold, it follows that
$\delta_{x}\per \theta -\psi_{*} \delta_{x}\per \theta=0$.
Since $\delta_{x}$ is arbitrary, it follows that $\psi^{*} \theta =\theta$.
That is, $\psi$ is a contact diffeomorphism.
 The converse is proved similarly. \qed

\begin{thm}
\label{thm:main}
Suppose that the Poisson manifold $P$ is compact.
Then, the Lie algebra of $\eex$ is the Poisson  algebra  $C^{\infty}(P)$.
\end{thm}
\pf  The Lie algebroid of $E$ is  the  central
extension of the cotangent Lie algebroid $T^*P$ by the Poisson
tensor $\pi$,  considered as a Lie algebroid 2-cocycle \cite{Xu:1994cmp}.
More precisely, the Lie algebroid of $E$ is
isomorphic to  $\tilde{A}=T^{*}P\oplus (P\times \reals)$, with the bracket being given by:
$$[(\zeta , f), (\eta , g)]=([\zeta  , \eta ], (\pi^{\#}\zeta )g- 
(\pi^{\#} \eta )f +\pi ( \zeta , \eta )), $$
for  any $\zeta  ,\eta \in \Omega^{1}(P) $, and $f, g\in C^{\infty}(P)$.
Here  a  pair $(\zeta  , f)$ for $\zeta  \in  \Omega^{1}(P) $ and 
$f \in C^{\infty}(P)$, 
is identified with a right   invariant vector field in $E$  via
the map $\lambda$: 
$$(\zeta  , f) \lon \widehat{X_{\alpha_{0}^{*}\zeta }}+
(p^{*}\alpha_{0}^{*}f)\xi,$$
where $\xi$ is the Euler vector field on $E$ generating
the  circle action.
 Let $\calg$ denote the group of all bisections
of $E$. Then, the  Lie algebra of $\calg $ can be
 identified with the space of sections of $\tilde{A}$ 
according to a  general rule. 
More explicitly, the tangent space $T_{1}\calg$  at the unit is identified
with $\gm (\tilde{A})$ as follows.   Any  isotopy $K_{t}$ 
in $\calg$
starting from $1$ 
 induces a  family of diffeomorphisms  on $E$ by the 
left  translations  $l_{K_{t}}$. The derivative of such a family of
diffeomorphisms at $t=0$ defines a  left invariant vector field $X$
on $E$, which can be canonically identified with a section  $(\zeta , f)$
of $\tilde{A}$ by the map $\lambda$ above.
 
Now if $K_{t}$ is  an isotopy
 in $\eex$, the corresponding vector field $X$ will be
 contact according to Proposition
\ref{pro:contact}.
If $\lambda (\zeta , f) =X$, then $X=\widehat{X_{\alpha_{0}^{*}\zeta }}+
(p^{*}\alpha_{0}^{*}f)\xi$.  It is easy
 to see that $X$ is a contact vector field, if and
only if $\zeta=-df$.
Therefore, the Lie algebra of $\eex$ can be  identified with the Lie subalgebra
of $\gm (\tilde{A})$ consisting of all elements of
the form $(-df, f)$, which is clearly isomorphic to
the Poisson algebra $C^{\infty}(P)$ via the identification
$(-df, f)\lon -f$. This concludes our proof. \qed

The following result    is then immediate.

\begin{cor}
\begin{enumerate}
\item The projection  $p: E\lon \gm$ induces
a group homomorphism $\eex \lon \gmm$, denoted by
the same symbol $p$, which lifts the Lie
algebra homomorphism $d: C^{\infty}(P)\lon Z^{1}(P)$.
\item $Im p=U_{ex}(\gm )$.
\end{enumerate}
\end{cor}
\pf For any bisection  $K$ of the groupoid $E$, $p K$ is clearly  a
bisection of $\gm$. Since $p : E\lon \gm$ is  a
groupoid morphism \cite{WeinsteinX:1991}, it follows
that $p \smalcirc l_{K}=l_{p K}\smalcirc p$.
If $K$ is a legendre  bisection, we have $ (l_{K})^{*}\theta =\theta $
according to   Proposition \ref{pro:contact}.
By taking the exterior derivative, one
gets that  $(l_{K})^{*} p^{*} \omega =p^{*} \omega$.
It  thus follows that $p^{*} l_{p K}^{*}\omega =p^{*} \omega$,
which implies  that $l_{p K}^{*}\omega=\omega$ since
$p$ is a submersion.
Therefore, $p K$ is a lagrangian  bisection.
This proves that we indeed get a map 
$p : \eex \lon \gmm$. This is a group homomorphism 
since $p: E\lon \gm$ is a groupoid morphism.
The last part of the proof of Theorem \ref{thm:main}
has already indicated  that the derivative of 
$ p : \eex \lon \gmm$  is indeed the exterior differential.

To prove the second assertion, let $K$ be any bisection
in $\eex$. Then there is an isotopy $K_{t}$ in $\eex$
such that $K_{0}=1$ and $K_{1}=K$. Let $L_{t}=p K_{t}\in U_{0} (\gm )$.
Then, $\Phi (L_{t})=0$ since $K_{t}$ is horizontal.
That is, $L_{t}$ is in $Ker \Phi$. Therefore
$L_{1}\in (Ker \Phi )_{0}$. However, 
$(Ker \Phi )_{0} =U_{ex}(\gm )$ according to Corollary \ref{cor:component}.
This shows that $Im p\subseteq U_{ex}(\gm )$. 

Conversely, for any
 given $L \in U_{ex}(\gm )= (Ker \Phi )_{0}$,
we assume that $L_{t}$ is an isotopy in $Ker \Phi$ connecting $1 $
and $L$. Then we can always find an isotopy  of bisections
$K_{t}$ in $E$ which is a parallel lifting  of $L_{t}$ since $\Phi (L_{t})=0$.
Then $K_{t}$ is an isotopy of legendre bisections.
That is, $K_{t}\in \eex$, and in particular $K_{1}\in \eex$.
Hence, we obtain the other inclusion: $U_{ex}(\gm )\subseteq Im p$.
This concludes the proof. \qed

Now we are ready to state the main theorem of the section.
\begin{thm}
When $P$ is a compact Poisson manifold
and its symplectic groupoid $\gm $ is
prequantizable, we have the following  exact sequence of group homomorphisms:

$$0\lon T^{1}\lon \eex\stackrel{p }{\lon} U_{0}(\gm )
\stackrel{\tilde{\frakf}}{\lon }H^{1}(P, \reals )/\Omega  \lon 0.$$

On the Lie algebra level, this corresponds to the 
exact sequence:

$$0\lon \reals \lon C^{\infty}(P)\stackrel{d}{\lon}Z^{1}(P)
\stackrel{f}{\lon } H^{1}(P, \reals ) \lon  0.$$
\end{thm}
\pf It remains to prove that 
 $ker \tilde{ \frakf } =U_{ex} (\gm )$. 

We  already know   that $U_{ex} (\gm ) \subseteq    ker \tilde{ \frakf }$.
For the other direction, let us assume that $L\in ker \tilde{ \frakf }$.
Let $L_{t}$ be any  isotopy  connecting $1$ and $L$.
Then, by the definition of $\tilde{\frakf}$, there  exists
an isotopy $S_{t}$ in $U_{0}( \gm )$ with $S_{1}=S_{0}=1$
such that  $\frakf (\{L_{t}\})= \frakf (\{S_{t}\})$.
>From this, it follows that $\frakf (\{L_{t}S_{t}^{-1} \} )=0$.
Therefore, $L_{1}S_{1}^{-1}$ is an exact bisection, so is $L=L_{1}$. \qed
  {\bf Remark.}   It is worth pointing out that in general a prequantum
bundle can be constructed whenever  the  periodic group
of $\gm$ is discrete. All the discussion in this section   could  be
carried out similarly    in this  generalized context.

\section{The flux homomorphism 
on $\alpha$-simply connected symplectic groupoids}

In this section,  again we will  assume that $(P, \pi )$ is a compact Poisson
manifold. Suppose that 
  $\nu$ is  a finite measure on 
 $P$, which is   invariant under all hamiltonian
flows.  For instance,
if $P$ is a compact symplectic
manifold with symplectic form $\omega$,
one can take  the
volume form $\omega^n$ as such a measure.
  Then $\nu$ is invariant under all
 vector fields of the form
$X_{\theta}$,   $ \forall \theta \in Z^{1}(P)$. 
Therefore, it is also invariant under the adjoint
action $Ad_{L}: P\lon P$, for any  $L\in U (\gm )$. 
Define a pairing 
$H^{1}(P, \reals )\times H^{1}(P, \reals )\lon \reals$ by 
\begin{equation}
\label{eq:pairing}
<[\theta_{1}], [\theta_{2}] >=\int_{\nu }\pi (\theta_{1}, \theta_{2} ).
\end{equation}

\begin{pro}
This pairing is well defined, skew-symmetric and bilinear.
\end{pro}
\pf It is clear that the rhs of Equation (\ref{eq:pairing}) is
bilinear and skew-symmetric with respect to
the arguments $\theta_{1}$ and $\theta_{2}$.
To prove that its value only depends on their cohomology
classes, it suffices to show that 
it vanishes when  either $\theta_{1}$ or $\theta_{2}$
is exact.
Assume that $\theta_{1}=df$. 
$$<[\theta_{1}], [\theta_{2}] >=\int_{\nu }\pi (df , \theta_{2} )
=-\int_{ \nu } X_{\theta_{2}}(f).$$
Let $\phi_{t}$ be the   flow   on $P$ generated
by $X_{\theta_{2}}$. Since  $\nu$  is invariant under $\phi_{t}$,
we have $\int_{\nu }\phi_{t}^{*} f=\int_{\nu }f$. 
By taking the  derivative, one gets that  $\int_{\nu } X_{\theta_{2}}(f)=0$.
This concludes the proof. \qed

 The pairing   $<\cdot , \cdot >$ induces a linear map
 $\mu: H^{1}(P, \reals )\lon
H^{1}(P, \reals )^*$.  Consider its  composition with the flux
homomorphism: $\varphi=\mu \smalcirc \frakf :\widetilde{\gmm} \lon H^{1}(P, \reals )^*$.
In general, the flux homomorphism $\frakf$ is  very   difficult
to compute. However,  the map $\varphi$ possesses  a nice description
 when $\gm$ is $\alpha$-simply connected (i.e. the $\alpha$-fibers
of $\gm$ are all simply-connected), as we will see below.
 This   description will   allow 
us to extract certain  information about the flux
homomorphism $\frakf$ even when an  explicit
expression  is  not available.  This is
particularly useful  when $\mu$ is an isomorphism. In this  case,  one  can
actually compute the flux homomorphism by first computing the map
$\varphi$.

To proceed, first  we  introduce a map
$\epsilon: U(\gm )\lon H^{1}(\gm  ,\reals )^*$  as follows:
$$<\epsilon (L), [J ] >=\int_{\nu } J(L(u)), \ \ J\in Z^{1}(\gm , \reals ), $$
where $H^{1}(\gm  ,\reals )$ is the first groupoid cohomology 
group  with real coefficients (see \cite{WeinsteinX:1991}).
 
\begin{pro}
\begin{enumerate}
\item The map $\epsilon$ is a well-defined map.
\item $\epsilon $ is a group homomorphism.
\end{enumerate}
\end{pro}
\pf Suppose that $J$ is a groupoid coboundary,
i.e., $J(x)=H(\alpha (x))-H(\beta (x))$ for some 
$H\in C^{\infty}(P)$.
Then, $J(L(u))=H(Ad_{L}u)-H(u)$ and
$\int_{\nu}J(L(u))= \int_{Ad_{L}^{*}\nu} H- 
\int_{\nu}H=0$, since
the measure $\nu$ is invariant. It is also simple to see
that $<\epsilon (L), J>$ is linear with respect  to $J$. Thus, $\epsilon $
is indeed  a well-defined map from $U(\gm )$ to   $H^{1}(\gm  ,\reals )^*$.

To prove the second assertion,
let $L$ and $K$ be  any two   lagrangian bisections.
Then, $\forall J\in Z^{1}(\gm , \reals  )$,
\be
<\epsilon (LK), [J]>&=&\int_{\nu}J(L(Ad_{K}u)K(u))\\
&=&\int_{\nu}J(L(Ad_{K}u))+\int_{\nu}J(K(u))\\
&=&\int_{Ad_{K}^{*}\nu}J(L(u))+\int_{\nu}J(K(u))\\
&=&\int_{\nu}J(L(u))+\int_{\nu}J(K(u))\\
&=&<\epsilon (L), [J]>+<\epsilon (K), [J] >. 
\ee
This concludes the proof. \qed

According to Theorem 1.3 in \cite{WeinsteinX:1991},
$H^{1}(\gm , \reals )$ is isomorphic to the first Poisson cohomology
$H^{1}_{\pi }(P)$ when $\gm$ is $\alpha$-simply
connected. On the cochain level, this isomorphism  $\Psi: H^{1}_{\pi }(P) \lon H^{1}(\gm , \reals )$
is established  by:
$$\Psi ([X] )(r)= \int_{\alpha (r)}^{r}\theta_{X}, \ \   \forall r\in \gm , $$
where $X$ is any  Poisson vector field on $P$, $\theta_{X}$ is its 
corresponding left-invariant  one-form on $\alpha$-fibers defined by
$\theta_{X} (X_{\beta^{*}f })=X(f), \ \forall f \in C^{\infty}(P)$,
  and the integration is
over any path   connecting  $\alpha (r)$ and $r$
in the $\alpha$-fiber.
On the other hand, the Poisson tensor  $\pi$ induces
a natural morphism $\pi^{\#}: H^{1}(P, \reals )\lon H^{1}_{\pi }(P)$.
By taking their duals and composing with $\epsilon $, we 
obtain a
map $\lambda : U(\gm )\lon H^{1}(P, \reals )^{*}$
such that  $\lambda =(\pi^{\#})^{*}
\smalcirc \Psi^{*} \smalcirc \epsilon $.

The following proposition gives an explicit  expression for this map
$\lambda$.

\begin{pro}
For  any $L\in U(\gm )$ and  $[\theta ]\in H^{1}(P, \reals )$,
$$<\lambda (L), [\theta ]>=\int_{\nu} ( \int^{\alpha^{-1}(u)\cap L}_{u} \beta^{*}
\theta ) , $$
where the interior integration is over any path in the $\alpha$-fiber
connecting the point $u$ and  its intersection with $L$:
the point $\alpha^{-1}(u)\cap L$.

\end{pro}
\pf Given any  closed one-form $\theta$, let $X=\pi^{\#}\theta $ and
$J=(\Psi \smalcirc \pi^{\#} )\theta $.

Since
$$(\beta^{*}\theta )(X_{\beta^{*}f })=\theta (T\beta X_{\beta^{*}f })=
-\theta (X_{f})= (\pi^{\#}\theta )f= Xf, $$
the restriction of $\beta^{*}\theta $ to $\alpha$-fibers
 will  be the left-invariant one-form  corresponding
to the Poisson vector field $X$.
Therefore, $\forall \gamma \in \gm$,

$$J(\gamma )=((\Psi \smalcirc \pi^{\#} ) \theta )(\gamma )=\Psi ([X] )(\gamma )
=\int_{\alpha (\gamma )}^{\gamma } \beta^{*}\theta .$$
Hence,
\be
<\lambda (L), [\theta ]>&=&<\epsilon (L) , (\Psi \smalcirc \pi^{\#} )[\theta ]>
\\
&=&<\epsilon (L) , [J]>\\
&=&\int_{\nu}J(L(u))\\
&=&\int_{\nu} (\int_{Ad_{L}u}^{L(u) } \beta^{*}\theta  )\\
&=&\int_{\nu} ( \int^{\alpha^{-1}(u)\cap L}_{u} \beta^{*} \theta )  ,
\ee
where the last step follows from the invariance of the measure $\nu$. \qed

We can now describe the main theorem of this section.
\begin{thm}
\label{thm:phi-lambda}
For any isotopy  $L_{t}$ in $\gmm$ connecting $1$ and $L$, we have
$$\phi (\{L_{t}\})=\lambda (L).$$
\end{thm}
\pf Given  any $[\xi ]\in H^{1}(P, \reals )$,  let 
$f(t)=<\lambda (L_{t}), [\xi ]>=\int_{\nu}J(L_{t}(u))$,
where $J=\Psi  (\pi^{\#}\xi )$ is the corresponding  groupoid 1-cocycle.
By taking the  derivative, one has
\be
\dot{f}(t)&=&\int_{\nu}X_{t}(J)(L_{t}(u)) \\
&=&-\int_{\nu}<X_{J}, X_{t}\per \omega >(L_{t}(u))\\
&=&-\int_{\nu}<X_{J}, \alpha^{*} \theta_{t}>(L_{t}(u))\\
&=&-\int_{\nu}<T\alpha X_{J}, \theta_{t}>(Ad_{L_{t} } u)\\
&=&-\int_{\nu}<\pi^{\#}\xi , \theta_{t}>(Ad_{L_{t}}  u)\\
&=&\int_{\nu}\pi (\theta_{t}, \xi )(u),
\ee
where the   second from the last  equality   follows from  
the proof of Proposition 2.2 in \cite{WeinsteinX:1991}.

Therefore,
$$f(1)=f(1)- f(0)=\int_{0}^{1}\dot{f}(t) =\int_{0}^{1}\int_{\nu}\pi (\theta_{t}, \xi )(u)dt.$$
On the other hand, it follows from the definition of $\phi$
that
\be
<\phi (\{L_{t}\}), [\xi]> &=&< \mu ( \frakf \{L_{t}\}), [\xi ]>\\
&=&\int_{\nu}\pi (\frakf (\{L_{t}\}), \xi ) \\
&=& \int_{0}^{1}\int_{\nu}\pi (\theta_{t},
 \xi )(u)dt .
\ee
This concludes the proof of the theorem. \qed

 An immediate consequence  of Theorem \ref{thm:phi-lambda} is the following:

\begin{cor}
\label{pro:descend}
When $\gm$ is $\alpha$-simply connected,
the map $\phi$ descends to  a group homomorphism
$ \tilde{\phi}: \gmm \lon H^{1}(P, \reals )^*$.
That is, $\phi$ vanishes on   $\pi_{1}(\gmm , 1)$.
\end{cor}

Another consequence is:

\begin{cor}
Suppose that $\gm $ is $\alpha$-simply connected.
For any $\{L_{t}\}\in \widetilde{U_{0}(\gm )}$ and $[\theta ]
 \in H^{1}(P, \reals )$,
$$<\frakf (\{L_{t}\} ), \mu [\theta ]>=- \int_{\nu} ( \int^{\alpha^{-1}(u)\cap L_{1}}_{u} \beta^{*}
\theta  ).$$
\end{cor}

When $\mu $ is an isomorphism, the above corollary  will  enable  us 
to carry out an explicit computation of the flux homomorphism.

Let us conclude this section by the following

\begin{cor}
\begin{enumerate}
\item $\gmex \subseteq Ker \lambda $.
\item In particular, if $\mu$ is an isomorphism,  then
$$\gmex=Ker  ( \lambda |_{ \gmm} ) . $$
\end{enumerate}
\end{cor}


\section{The case of symplectic manifolds}
As an example, we will consider compact symplectic manifolds in this
section. In particular, we will investigate
the relation between the  flux homomorphism of  their
symplectic groupoids (being the fundamental groupoids in this case)
and the usual flux homomorphism of symplectomorphisms.

Let $M$ be a  compact symplectic manifold. Then, $M$ is equipped with a
natural invariant volume form $\omega^n$.
We can take $\nu$ to be its corresponding 
measure.
The symplectic pairing,
as  often called in the literature  \cite{Banyaga}  \cite{Calabi} \cite{McduffS}, is the pairing
in $H^{1}(M, \reals )$ given by: 
\begin{equation}
\sigma ([\theta_{1}], [\theta_{2}] )=\int \theta_{1}\wedge \theta_{2}\wedge
\omega^{n-1},   \ \ \forall \theta_{1} , \theta_{2}\in H^{1}(M, \reals ).
\end{equation}

It is well known that the symplectic pairing is nondegenerate
in the case of the symplectic structure subordinate to a
compact Kahler manifold.
The following lemma can be proved easily by using a  local
Darboux chart.

\begin{lem}
Up to a constant, $\sigma$ coincides with the pairing $<\cdot , \cdot >$
defined by Equation (\ref{eq:pairing}).
\end{lem}

For $M$, its   $\alpha$-simply connected 
symplectic groupoid  $\gm $ is the fundamental groupoid $\Pi_{1}(M)$. 
For any $[\theta ]\in H^{1} (M, \reals )$, let 
$J\in C^{\infty}(\gm )$ be given by,
$$ J(\gamma )=\int_{\gamma }\theta  ,  \ \ \ \ \forall \gamma \in \gm =\Pi_{1} (M).$$
Then
$(\Psi \smalcirc \pi^{\#}
)[\theta ]=[J]$ according to   \cite{MikamiW:1988} (in this case, 
   $(\Psi \smalcirc \pi^{\#}) $ is an isomorphism).
Therefore, $\lambda : U_{0} (\gm )\lon H^{1} ( M  , \reals  )^{*}$ is
given by
$$<\lambda  (L), [\theta ]>=\int_{M} (\int_{L(u)}\theta )\omega^{n}
=\int_{[0, \ 1]\times M}i^{*}(\theta_{x}\wedge \omega^{n}_{y}),$$
where $i$ is the embedding $[0, \ 1]\times M\lon M\times M$ defined as
$(t, u)\lon (L(u)(t), u)$,  and $\theta_{x}$ is  the
one-form $\theta $ considered as a one-form
on the first component while $\omega_{y}$ is $\omega $  considered as  a two-form
 on the second component.

Given any $L\in U_{0} (\gm )$.  For any fixed $u\in M$, $L(u)$ is an element
in $\Pi_{1} (M)$ with the end  point $u$.
By an isotropic bisection, we mean a bisection of $\gm$  whose all points
 lie in the isotropic groupoid of $\gm$.
 Given any  isotropic bisection $L$ in $U_{0} (\gm )$, $L(u)$
is a homotopy class of a closed path  with the  end point $u$.
It is simple to see that for any $u$ and $v$ in $M$, $L(u)$
and $L(v)$ are homotopy equivalent loops  since $M$ is
connected.  Therefore, $\int_{L(u)}\theta =\int_{L(v)}\theta $, 
or $J(L(u))=J(L(v ))$. It thus follows that
$$<\lambda (L), [\theta ]>=\int_{M}J(L(u))\omega^{n}=J(L (u)) \cdot 
Vol (M)=
(\int_{L(u)}\theta  )\cdot  Vol(M). $$

An immediate consequence of this formula is the following:

\begin{pro}
\label{pro:null}
An isotropic lagrangian bisection $L$ is in $ker \lambda$
iff $L(u)$ is null-homologous for every $u\in M$.
In particular, if  an isotropic bisection $L\in U_{0} (\gm )$ is exact, then $L(u)$ is null-homologous.
\end{pro}

Given any symplectomorphism $f$, by the graph of $f$, we mean
the submanifold $\{(f(x), x)|x \in M\}$ of $M\times \Bar{M}$.
Then the  graph of a symplectomorphism is clearly a lagrangian
submanifold.
Let $p$ denote the  projection $p=\alpha \times \beta :
\Pi_{1}(M)\lon M\times M$.
The following result reveals the relation between the
usual flux homomorphism of symplectomorphisms
and the flux homomorphism of symplectic groupoids.

\begin{pro}
Let $h_{t}$ be  a symplectic isotopy  in $Diff_{0}(M, \omega )$
 such that $h_{0}=id$.
Suppose that $h_{t}$ lifts to an isotopy of   bisections  $L_{t}$
in  $\gm$ with $L_{0}=1$  (i.e.,  $L_{t}$ goes to the graph 
of $h_{t}$ under the projection $p$). Then,
\begin{enumerate}
\item $L_{t}$ is a  lagrangian bisection;
\item $ \frakf (\{L_{t}\})= F(\{h_{t}\})$.
\end{enumerate}
\end{pro}
\pf Let $\tilde{\omega }$ be the symplectic form on $\gm$. Then,
$\tilde{\omega }=\alpha^{*}\omega -\beta^{*}\omega $. Let $i$ and
$i'$ denote  the inclusions: $i: L_{t}\lon \gm $ and $i': 
\mbox{ (the graph of }\  h_{t})\lon M
\times \Bar{M}$,  respectively. Then,
$i^{*}\tilde{\omega }=(p|_{L_{t}})^{*} (i')^{*} (\omega \ominus \omega )=0$,
since the  $ \mbox{graph of }  h_{t}$ is  lagrangian. This shows that $L_{t}$ is
a lagrangian bisection.

Let $Y_{t}$ be the  corresponding time-dependent  vector field on $\gm $ 
defined by Equation (\ref{eq:Y}), and $\theta_{t}$
the    time-dependent one-form on $P$ as  defined   by Equation (\ref{eq:one-form}).
Thus,
\be
Y_{t}\per \tilde{\omega }&=&(Y_{t}\per \alpha^{*}\omega )-
(Y_{t}\per \beta^{*}\omega )\\
&=&\alpha^{*}(T\alpha Y_{t} \per \omega ).
\ee
Since $\alpha L_{t}(m)=h_{t}(m)$ for any $m\in M$, 
it follows that $T\alpha Y_{t}=\dot{h_{t}}$.
Hence, $Y_{t}\per \tilde{\omega }=\alpha^{*}(\dot{h_{t}}\per \omega )$.
In other words, $\theta_{t}=\dot{h_{t}}\per \omega$.
Thus, $\frakf (\{L_{t}\})=[\int_{0}^{1}\theta_{t}dt ]=[\int_{0}^{1}(\dot{h_{t}}\per \omega )dt ]
= F (\{h_{t}\})$. \qed

By combining  with Theorem  \ref{thm:phi-lambda}, this proposition leads 
to the following:
\begin{thm}
\label{thm:flux-nondeg}
For   $\{ h_{t}\}\in \widetilde{  Diff_{0}(M, \omega ) }$,
suppose that $L_{t} $ is any  isotopy in $\gmm$ lifting
$\{h_{t}\}$.
Then, 
\begin{equation}
\label{eq:mu-F}
(\mu  \smalcirc F)\{h_{t} \} =\lambda (L_{1}). 
\end{equation}
In particular, if the symplectic pairing   is nondegenerate,  we have
$$F (\{ h_{t} \})=\mu^{-1}\lambda (L_{1} ). $$
\end{thm}

This theorem suggests a useful method for  computing 
the flux homomorphism for symplectomorphisms. 
The rhs of Equation (\ref{eq:mu-F}) can be essentially
expressed by an  integration over the symplectic
manifold, which is much easier to handle. 
In the case when $\mu $ is an isomorphism, this
will enable us to carry out an explicit 
computation for  the usual flux homomorphism $F$.

We note that
 any     symplectic isotopy $h_{t}$ of  $M$  can always be lifted to an
isotopy of lagrangian bisections in $\gm$.
 For example,  for any fixed $t$,
set $L_{t}=\{[h_{t-s}(m)], 0\leq s\leq t| m\in M\}$.
Then it is simple to see that $L_{t}$ is a lifting of $h_{t}$.

Let $A: \widetilde{Diff_{0}(M, \omega )}\lon \gmm$ be the map  which
sends    
$\{h_{t}\} $ to $L=\{[h_{1-s}(m)],\ 0\leq s\leq 1| m\in M\}$.

The following corollary follows directly from Theorem 
\ref{thm:flux-nondeg}.

\begin{cor}
The following diagram:
\begin{equation}                         
\matrix{&& A &&\cr
        &  \widetilde{  Diff_{0}(M, \omega ) } &\vlra&  \gmm
 &\cr
        &&&&\cr
     F &\Bigg\downarrow&&\Bigg\downarrow&\lambda  \cr
        &&&&\cr
        & H^{1}(M, \reals ) &\vlra&  H^{1}(M , \reals )^{*} &\cr
        && \mu  &&\cr}
\end{equation}
commutes.

\end{cor}

Combining  Proposition \ref{pro:null} with
 the corollary above  leads  
 the  following result of McDuff \cite{Mcduff}:

\begin{cor}
A loop of symplectomorphisms $\{\phi_{t }\}_{ 0\leq t\leq 1} \in \pi_{1}
 (Diff_{0}(M, \omega ), 1 )$
  is in the
kernel of $\mu \smalcirc F$ if and only if  the loop
$u\lon \phi_{t}(u) \in M$ is null-homologous for every $u\in M$.
\end{cor}
{\bf Remark} (1).  Clearly,  $A$ maps the fundamental group
 $\pi_{1}(Diff_{0}(M, \omega ), 1)$ into
the group  $I$ of isotropy lagrangian bisections of $\gm$.
 It would
be interesting to know if $A$ is surjective. It
is also interesting to ask if $I$ is discrete.

(2). Suppose that $\phi_{t}$ with $\phi_{1}=\phi_{0}=id$  is
a loop in $Ham (M, \omega )$ for  a symplectic manifold $P$.
It was asked in \cite{McduffS} whether the loop $u\lon \phi_{t}(u)$ is
contractible for every $u\in P$.  We may ask
the following more general   question:  does  the  intersection of
$U_{ex}(\gm )$ with $I$ consist only the trivial element $1$.
Note that this fails for  general Poisson manifolds.
For instance,  for zero Poisson manifolds, it is simple to see that
$U_{ex}(\gm )\cap I= U_{ex}(\gm )$.

\begin{pro}
Suppose that the image of $\pi_{1}(Diff_{0}(M, \omega ), 1)$ under $A$
coincides with the space  $I$ of all isotropy lagrangian bisections and
$\mu $ is an isomorphism.
Then, the flux homomorphism $F$ descends to a homomorphism
$F': Diff_{0}(M, \omega )\lon H^{1}(M, \reals )/\Omega $,
where $\Omega =(\mu^{-1}\lambda )(I)$ and the map $F'$
is given by
$F' (h)=[(\mu^{-1}\lambda ) (L)]\in H^{1}(M, \reals )/\Omega  $,
$\forall h\in Diff_{0}(M, \omega )$.
Here $L$ is any  bisection in $\gmm$ such that  $p(L)=
\mbox{the graph of  }  h$.
\end{pro}

We end this section by applying the results above
to the symplectic torus $T^{2n}$,  which was treated
in \cite{McduffS} by a different method.

\begin{ex}

Consider the torus $T^{2n} \cong \reals^{2n}/\integers^{2n}$ with
the  canonical symplectic structure.
Let $\phi_{t}$ be a symplectic isotopy with a lift $\psi_{t}: \reals^{2n}
\lon \reals^{2n}$ such that $\psi_{t} (x+l)=\psi_{t}(x)+l, \ \forall
l\in \integers^{2n} $ and
$\psi_{0}=id$ for any $x\in \reals^{2n}$. The symplectic groupoid
$\gm$ can be taken as  the cotangent space $T^{*} T^{2n} \cong  T^{2n}\times
\reals^{2n}$ with the standard symplectic structure and
the source and target maps are given respectively
by $\alpha (x, p)=x-\half p$, $\beta (x, p)=x+\half p$ (see  \cite{Weinstein:1991sg}).
Clearly, the first groupoid cohomology $H^{1}(\gm  , \reals )$ is generated by
$J_{i}(x, p)=p_{i}, \ 1\leq i\leq 2n$. Take
$L_{t}=\{(\half (x+\psi_{t}(x)), x-\psi_{t}(x))\in  T^{2n}\times
\reals^{2n}|x \in \reals^{2n}\}$.
 It is simple to see that
both $\alpha$ and $\beta $ are diffeomorphisms when restricted  to $L_{t}$ and
$p(L_{t})=graph (\phi_{t})$.
In other words, $L_{t}$ is a family of lagrangian bisections of
$\gm$, which lifts the graph of $\phi_{t}$.
Now,
$$<\epsilon  (L_{1}) , [J_{i}] >=\int_{T^{n}}J_{i}(L_{1}(x))=
\int_{T^n}(x^{i}-\psi_{1}^{i}(x))dx .$$
On the other hand, it is simple to see that $(\Psi \smalcirc \pi^{\#})[dx_{i}]
=[J_{i}]$, $i=1, \cdots , 2n$.
Therefore, 
$$<\lambda  (L_{1}), [dx_{i}]>
=<\epsilon  (L_{1}), (\Psi \smalcirc \pi^{\#})[dx_{i}]>=
<\epsilon  (L_{1}), [J_{i}]>=
\int_{T^n}(x^{i}-\psi_{1}^{i}(x))dx  . $$
Hence, according to Theorem    \ref{thm:flux-nondeg}, 
the flux homomorphism $F$ maps $\{\phi_{t}\}$    to
$[\Sigma_{j=1}^{2n}a_{j}dx_{j}]$,
where $a=(a_{1}, \cdots ,a_{2n})= (\int_{T^{2n}}(x^{1}-\psi_{1}^{1}(x))dx, \cdots ,
\int_{T^{2n}}(x^{2n}-\psi_{1}^{2n}(x))dx) \cdot  J $, 
and
$J$ is the canonical symplectic matrix on $\reals^{2n}$.
\end{ex}
{\bf Remark} The $\alpha$-simply connected 
 symplectic groupoid  of  a compact 
  Kahler manifold  $P$ with  negative holomorphic constant curvature
is isomorphic to the cotangent bundle $T^{*}P$ \cite{Weinstein:pp}. It would be
interesting to generalize the computation
  above to
such a manifold.  It is reasonable to expect that
in this case  the flux homomorphism
should be related to certain
``center of mass''  of the  manifold.

     \end{document}